\documentclass[12pt]{amsart}

\usepackage{graphicx}
\usepackage{amsfonts}
\usepackage{amssymb}
\usepackage{hyperref}
\usepackage{url}

\makeatletter
\def\url@leostyle{%
  \@ifundefined{selectfont}{\def\UrlFont{\sf}}{%
    \def\UrlFont{\footnotesize\ttfamily}}\Url@do
} \makeatother

\newtheorem{lemma}{Lemma}
\newtheorem{theorem}[lemma]{Theorem}
\newtheorem{proposition}[lemma]{Proposition}
\newtheorem*{proposition*}{Proposition}
\newtheorem*{theorem*}{Theorem}

\theoremstyle{definition}
\newtheorem{definition}{Definition}

\newtheorem{conjecture}[definition]{Conjecture}

\theoremstyle{remark}

\newcommand{\ci}{\subseteq}

\newcommand{\intersect}{\cap}
\newcommand{\avoid}{{\mathrm{Av}}}
\newcommand{\juxt}{{\mathrm{Juxt}}}
\newcommand{\merge}{{\mathrm{Merge}}}
\newcommand{\cycle}{{\mathrm{Rot}}}

\newcommand{\expectation}{{\mathbf{E} \,}}
\newcommand{\probability}{{\mathbf{Pr} \,}}

\newcommand{\LS}[1]{{\mathrm{L} {#1} \mathrm{S}}}
\newcommand{\A}{{\mathcal A}}
\newcommand{\B}{{\mathcal B}}
\newcommand{\C}{{\mathcal C}}
\renewcommand{\L}{{\mathcal L}}
\newcommand{\U}{{\mathcal U}}
\newcommand{\J}{{\mathcal J}}
\newcommand{\M}{{\mathcal M}}
\newcommand{\R}{{\mathcal R}}

\newcommand{\I}{{\mathcal I}}
\renewcommand{\S}{{\mathcal S}}
\newcommand{\sw}{s}

\title[Longest pattern avoiding subsequences]{On the length of the longest subsequence avoiding an arbitrary pattern in a random permutation}
\author{M.~H.~Albert}
\address{Department of Computer Science \\%
University of Otago \\
Dunedin, New Zealand} \email{malbert@cs.otago.ac.nz}

\begin{document}

\begin{abstract}
We consider the distribution of the length of the longest subsequence avoiding an arbitrary pattern, $\pi$, in a random permutation of length $n$. The well-studied case of a longest increasing subsequence corresponds to $\pi = 21$. We show that there is some constant $c_{\pi}$ such that as $n \to \infty$ the mean value of this length is asymptotic to $2 \sqrt{c_{\pi} n}$ and that the distribution of the length is tightly concentrated around its mean. We observe some apparent connections between $c_{\pi}$ and the Stanley-Wilf limit of the class of permutations avoiding the pattern $\pi$.
\end{abstract}

\maketitle

\section{Introduction}

The aim of this paper is to generalize certain results concerning the longest increasing subsequences of permutations to longest subsequences avoiding some pattern or patterns. The former area is quite well known, while the latter, though an active field of current combinatorial research, is considerably less so. Thus before surveying the known results (or at least the ones relevant to this paper) concerning longest increasing subsequences, we introduce the fundamental definitions concerning pattern avoidance in permutations.

Let $\sigma = \sigma_1 \sigma_2 \cdots \sigma_k$ and $\pi = \pi_1 \pi_2 \cdots \pi_n$ be permutations written as sequences. Then $\sigma$ {\em occurs as a pattern\/} in $\pi$ (or $\sigma$ is {\em involved in \/} $\pi$) if some length $k$ subsequence of $\pi$ has the same relative ordering as $\sigma$. In other words, for some $1 \leq i_1 < i_2 \cdots < i_k \leq n$ and all $1 \leq s,t \leq k$:
\[
\sigma_s < \sigma_t \: \iff \: \pi_{i_s} < \pi_{i_t}.
\]
For example, a permutation  $\pi$ contains an increasing subsequence of length $k$ if and only if $123 \cdots k$ occurs as a pattern in $\pi$. We say that $\pi$ {\em avoids\/} $\sigma$ if $\sigma$ is not involved in $\pi$. A {\em pattern avoidance class\/}, $\avoid(B)$ is the set of all permutations that avoid each of the permutations in the set $B$. The set $\avoid(B)$ is infinite unless $B$ contains both an increasing and a decreasing permutation (in which case it is finite by the Erd\"{o}s-Szekeres theorem) and if $B$ is non empty then $\avoid(B) \intersect \S_n$ is a proper subset of $\S_n$ for all $n$ greater than or equal to the length of the shortest element of $B$.

More generally, given a finite sequence $\sigma$ of distinct values from any linearly ordered set, we define its {\em pattern\/} to be that permutation whose elements are in the same relative order as those of $\sigma$. Then involvement can be defined for such sequences by reference to their patterns (or directly as above).

Recently Marcus and Tardos \cite{MT:JCTA04} proved the Stanley-Wilf conjecture about pattern avoidance classes. Specifically,
\begin{theorem*}[Marcus-Tardos Stanley-Wilf] For any proper pattern avoidance class $\A$ there is a real number $\sw_{\A}$ such that:
\[
\limsup_{n \to \infty} \left| \S_n \intersect \A \right|^{1/n} = \sw_{\A}.
\]
\end{theorem*}
So, not only is $\avoid(B) \intersect \S_n$ a proper subset of $\S_n$ for sufficiently large $n$, it is in fact small relative to $\S_n$, being of merely exponential size in $n$. We call $\sw_{\A}$ the {\em Stanley-Wilf limit} of the class $\A$.

Consider an arbitrary pattern avoidance class $\A$. Given a permutation $\pi$ define the {\em longest $\A$ subsequences\/} of $\pi$ or $\LS{\A}(\pi)$ to be the set of those subsequences of $\pi$ of maximum length, subject to the condition that their patterns belong to $\A$. Also define $L_{\A}(\pi)$ to be the length of any sequence in $\LS{\A}(\pi)$. Let $\I = \avoid(21)$ be the class of increasing permutations (note that we should really write $\avoid(\{ 21 \})$ but we will frequently omit such braces for explicit sets of avoided patterns). Then $\LS{\I}(\pi)$ is simply the set of longest increasing subsequences of $\pi$.

Apparently Ulam \cite{Ulam:MC61} was the first to ask the question:
\begin{quote}
{\em What can be said about the distribution of values of $L_{\I}(\Pi_n)$ when $\Pi_n$ is a random variable whose value is a permutation $\pi$ chosen uniformly at random from among the elements of $\S_n$?}
\end{quote}

We intend to address the generalization of this problem to the random variable $L_{\A}(\Pi_n)$ defined in a similar fashion. The history of the analysis of Ulam's problem is well documented in \cite{AD:BAMS99}. We repeat here a few details relevant to our investigations of the more general problem.

For convenience let $L_n = L_{\I}(\Pi_n)$. Ulam conjectured that for some constant $c$
\[
\lim_{n \to \infty} \frac{\expectation L_n}{\sqrt{n}} = c.
\]

This conjecture was proven by Hammersley \cite{Hammersley:Seed72} who showed also that $n^{-1/2} L_n \to c$ in probability, and who conjectured that $c = 2$. This further conjecture was proven in part by Logan and Shepp \cite{LS:AM77} and simultaneously in whole by Kerov and Ver{\v{s}}ik \cite{KV:DAN77}. 

Frieze \cite{Frieze:AAP91} and Bollob{\'a}s and Brightwell \cite{BB:AAP92} used martingale methods to establish tight concentration of $L_n$ about its mean. Subsequently, Baik, Deift and Johansson \cite{BDJ:JAMS99} obtained complete asymptotic information about the distribution of $L_n$. 

The evaluation of $L_{\I}(\pi)$ for a specific permutation $\pi$ is straightforward, owing to the fact that it is the length of the first row of the Young tableau obtained from $\pi$ using the Schensted algorithm \cite{Schensted:CJM61}. For an arbitrary permutation $\pi$ this row is easily constructed in worst case time $O(|\pi| \log L_{\I}(\pi) )$. Furthermore, the Robinson-Schensted-Knuth correspondence allows exact computation of the distribution of $L_n$ for moderate values of $n$ as was done by Baer and Brock \cite{BB:MC68} and to larger values by Odlyzko and Rains \cite{OR:CM00}.   Monte Carlo simulation of the distribution of $L_n$ for much larger values of $n$ carried out in the early 1990's but reported in \cite{OR:CM00} led to new conjectures concerning the asymptotic variance of the distribution, conjectures subsequently proved correct in \cite{BDJ:JAMS99}. 

Our main theorems are analogs of the results of Hammersley and Frieze for the more general case of $L_{\A}(\Pi_n)$. For the first of these we need to impose a mild additional restriction on $\A$. The proofs for the general case are then essentially identical to the originals. We will illustrate an apparent connection between the limit of the expectation of the longest $\A$ subsequence and the constant $\sw_{\A}$ together with some simulation results for classes where we have been unable to compute the limiting ratio exactly.

\section{Longest $\A$ subsequences in random permutations}

Let a proper pattern avoidance class $\A = \avoid(B)$ be given. Choose some constant $s$ such that for all sufficiently large $n$, 
\[
\left| \A \intersect \S_n \right| < s^n.
\]
The existence of such a constant follows from the Marcus-Tardos Stanley-Wilf Theorem. 

As in the arguments of \cite{Frieze:AAP91} it is convenient to work with a different random variable than $L_{\A}(\Pi_n)$. Let $\mathbf{X} = (X_1, X_2, \ldots, X_n)$ be a sequence of independent uniform $[0,1]$ random variables. Define the random variable $L_{\A}(\mathbf{X})$ to be the length of the longest subsequence of $\mathbf{X}$ whose pattern belongs to $\A$. This random variable depends only on the pattern of the sequence $\mathbf{X}$ and these patterns are uniformly distributed over $\S_n$. 

The proof in  \cite{Frieze:AAP91} begins with a crude probability inequality for $L_n$. We prove essentially the same inequality in exactly the same way, allowing only for the growth rate of the class $\A$ (bounded by $s^n$) as opposed to that of $\I$ (bounded by 1).

\begin{lemma}
\label{LEM:Probability}
\[
\probability \left( L_{\A}(\mathbf{X}) \geq 2e \sqrt{sn} \right) < e^{-2e\sqrt{sn}}.
\]
\end{lemma}

\begin{proof}
Let $n_0 = \lceil 2e \sqrt{sn} \rceil$, and let $\sigma$ denote the number of subsequences of $\mathbf{X}$ of length $n_0$ whose pattern belongs to $\A$. Then:
\begin{eqnarray*}
\probability \left( L_{\A}(\mathbf{X}) \geq n_0 \right) &\leq& \expectation (\sigma) \\
{} &\leq& \frac{{n \choose n_0} s^{n_0}}{n_0!} \\
{} &\leq& \left( \frac{n e^2 s}{n_o^2} \right) ^{n_0} \\
{} &<& e^{-2e \sqrt{sn}}.
\end{eqnarray*}
\end{proof}

Define the direct sum $\alpha \oplus \beta$ of two permutations $\alpha$ and $\beta$ to be that permutation $\alpha \beta'$ where the pattern of $\beta'$ is the same as that of $\beta$ and each element of $\beta'$ is larger than all the elements of $\alpha$. Similarly the direct difference $\alpha \ominus \beta$ is defined to be that permutation $\alpha' \beta$ where the pattern of $\alpha'$ is the same as that of $\alpha$ and each element of $\alpha'$ is larger than all the elements of $\beta$. A pattern class $\A$ is said to be {\em sum-closed\/} (respectively {\em difference-closed\/}) if for all $\alpha, \beta \in \A$, $\alpha \oplus \beta \in \A$ (respectively $\alpha \ominus \beta \in \A$).

Much of the investigation of pattern avoidance classes has focused on the {\em principal} classes, those of the form $\avoid(\pi)$ for some single permutation $\pi$ of length at least two. Any principal pattern avoidance class is either sum or difference closed since if $\A = \avoid(\pi)$ then $\A$ is sum-closed (difference-closed) if $\pi$ has no representation of the form $\theta \oplus \tau$ ($\theta \ominus \tau$) with both $\theta$ and $\tau$ non-empty. However, no permutation can be written both as a proper sum and as a proper difference.

\begin{theorem}
\label{THM:Expectation}
Let $\A$ be an infinite and proper pattern avoidance class which is either sum-closed or difference-closed. There exists a constant $1 \leq c_{\A} < \infty$ such that
\begin{equation}
\label{eq01}
\lim_{n \to \infty} \frac{\expectation L_{\A}(\Pi_n)}{\sqrt{n}} = 2\sqrt{c_{\A}}.
\end{equation}
\end{theorem}

\begin{proof}
Without loss of generality assume that $\A$ is sum-closed (the set of reversals of permutations in a difference closed class is sum-closed). The proof is a direct modification of Hammersley's for $L_n$. We repeat a sketch of the details here.  A more complete exposition can be found in \cite{Durrett:Probability}. Take a Poisson point process in $(0,\infty) \times (0, \infty)$ with unit density per unit area. Let $A_k$ be the length of the longest $\A$ subsequence in the $k \times k$ box with lower left corner at $(0,0)$. Then because $\A$ is $\oplus$-closed:
\[
A_{k+m} \geq A_k + A_m \circ T^k
\]
where $T^k$ shifts the origin to $(k,k)$ and preserves measure. The subadditive ergodic theorem applies to $-A_k$ and implies that $A_k/k$ converges almost surely to some value in $[0,\infty]$. Conditioned on the number of points in the $k \times k$ square being $n$, $A_k$ is distributed as $L_{\A}(\Pi_n)$. For large $k$, $n$ is approximately $k^2$. Every sum-closed class contains all increasing permutations which establishes the lower bound claimed. Furthermore,  one consequence of the lemma above is that $\limsup_{n \to \infty} L_{\A}(\mathbf{X})/\sqrt{n} < \infty$ so the full result follows.
\end{proof}

Reasons for the obscure choice of limiting constant in equation \ref{eq01} will be given in the next section. We now present a modification of the arguments in \cite{Frieze:AAP91} to show that the expectation of $L_{\A}(\Pi_n)$ is tightly concentrated around its mean.

\begin{theorem}
\label{THM:Concentration}
Let $\A$ be a proper pattern class.
For $\alpha > 1/3$ and $\beta < \min ( \alpha, 3 \alpha - 1 )$
\begin{equation}
\label{eq02}
\probability \left( \left| L_{\A}(\Pi_n) - \expectation L_{\A}(\Pi_n) \right| \geq n^\alpha \right) < \exp (-n^{\beta}).
\end{equation}
\end{theorem}

\begin{proof}
With the probability bound supplied by the lemma above in place 
the remainder of the proof for tight concentration is virtually identical to that of  \cite{Frieze:AAP91}. At some points in that argument various constants appear (for example, $6$ whose significance is that it dominates $2e$). Owing to the difference in the probability bounds, these constants must all be multiplied by $\sqrt{s}$. Otherwise the argument is completely unaffected, since it never actually makes use of any properties of increasing sequences.
\end{proof}

We remark at this point that the reason we have not been able to give the tighter concentration bounds supplied by \cite{BB:AAP92} is that the arguments of that paper do use special properties of increasing sequences. Specifically, they make use of a decomposition of the unit square into subsquares and the fact that if an increasing subsequence contains an element from a particular subsquare, then this generally rules out elements in many other subsquares. However, the form of (\ref{eq02}) is borrowed from \cite{BB:AAP92} where the observation is made that this form follows directly from the proof provided in \cite{Frieze:AAP91}.

\section{Meaning of the constant}
\label{SEC:Theory}

The Marcus-Tardos Stanley-Wilf theorem was a key ingredient in the proof of Theorems \ref{THM:Expectation} and \ref{THM:Concentration}. We will now present observations which provide some evidence for a connection between the constants $c_{\A}$ and $\sw_{\A}$.

Recall the definition of $\sw_{\A}$:
\[
\sw_{\A} = \limsup_{n \to \infty} \left| \S_n \intersect \A \right|^{1/n}.
\]
A superadditivity argument due to Arratia \cite{Arratia:EJC99} establishes 
that the right hand side is an actual limit for any principal pattern avoidance class, and more generally for any sum or difference-closed pattern avoidance class. 

For convenience we adopt a similar definition for $c_{\A}$ that applies even if we cannot be sure that the limit given in equation \ref{eq01} exists.
\[
2 \sqrt{c_{\A}} = \limsup_{n \to \infty} \frac{\expectation L_{\A}(\Pi_n)}{\sqrt{n}}.
\]
Certainly Lemma \ref{LEM:Probability} implies that for an infinite and proper pattern class $\A$, $1 \leq c_{\A} < 6 \sw_{\A}$. 

\begin{conjecture}
\label{conj01}
For any proper pattern avoidance class $\A$, the limits superior definining $c_{\A}$ and $\sw_{\A}$ are in fact limits, and $c_{\A} = \sw_{\A}$.
\end{conjecture}

This conjecture provides the reason for the otherwise obscure choice for the form of the right hand side of equation \ref{eq01}. The similarity of the arguments required to prove that both limits superior above are actually limits in the sum or difference-closed case is our second piece of (weak) evidence for the conjecture (the first is that it applies to $\I$!) The remainder of the paper is devoted to an investigation of the rather fragmentary evidence in support of Conjecture \ref{conj01}. This evidence takes two forms: a series of propositions that establish the equality claimed by it holds for certain classes and is preserved by certain constructions, and some limited experimental data on two specific classes. We begin with the more theoretical evidence. 

It will be convenient to use the expression ``$\A$ satisfies Conjecture \ref{conj01}'' as an abbreviation for the more accurate ``the limits superior defining $c_{\A}$ and $\sw_{\A}$ are in fact limits, and $c_{\A} = \sw_{\A}$''. We will also consider occasionally the weak form of the above where only the equality of $c_{\A}$ and $\sw_{\A}$ is asserted and not the fact that these constants can be computed as limits rather than as limits superior.

\begin{proposition}
For each positive integer $k$, the classes 
\[
\avoid (k(k-1)(k-2) \cdots 21) \quad \mbox{and} \quad \avoid (123 \cdots k)
\]
satisfy Conjecture \ref{conj01}.
\end{proposition}

\begin{proof}
By symmetry it suffices to consider $\I(k) = \avoid (k(k-1)(k-2) \cdots 21)$. Greene \cite{Greene:AM74} proved that the length of the longest $\I(k)$-subsequence of an arbitrary permutation $\pi$ is the sum of the lengths of the first $k-1$ rows of the tableaux produced from $\pi$ through the Robinson-Schensted-Knuth algorithm. From the detailed asymptotics for the shapes of such tableaux as can be found in  \cite{BOO:JAMS00} we deduce that
\[
\lim_{n \to \infty} \frac{\expectation L_{\I(k)} (\Pi_n)}{\sqrt{n}} = 2(k-1)
\]
and therefore $c_{\I(k)} = (k-1)^2$. On the other hand results of Regev \cite{Regev:AM81} (and see also \cite{GWW:EJC98}) demonstrate that $\sw_{\I(k)} = (k-1)^2$.
\end{proof}

We now consider four constructions which when applied to pattern avoidance classes $\A$ and $\B$ define a pattern avoidance class $\C \supseteq \A \cup \B$. The first of these is the union construction itself, the second is direct sum (which we have seen above), the third is called {\em juxtaposition}, and the fourth {\em merge}.

The juxtaposition $\juxt (\A, \B)$ of two pattern avoidance classes consists of the set of all permutations of the form $\alpha \beta$ where the pattern of $\alpha$ lies in $\A$ and that of $\beta$ lies in $\B$. The merge $\merge (\A, \B)$ consists of the set of all permutations $\pi_1 \pi_2 \cdots \pi_n$ such that for some subset $I \ci \{1,2,\ldots, n\}$ the pattern of the subpermutation of $\pi$ consisting of elements whose indices lie in $I$ belongs to $\A$, and the pattern of the remaining elements lies in $\B$.

\begin{proposition}
Let $\A$ and $\B$ be two pattern avoidance classes which satisfy Conjecture \ref{conj01}. Then their union, direct sum and juxtaposition also satisfy Conjecture \ref{conj01}. If, additionally, $\A \cap \B$ is a finite class then their merge also satisfies Conjecture \ref{conj01}. 
\end{proposition}

\begin{proof}
Suppose that $\A$ and $\B$ are pattern avoidance classes which satisfy Conjecture \ref{conj01} and let $a_k$ and $b_k$ denote the number of permutations of length $k$ in $\A$ and $\B$ respectively. 

Let $\U = \A \cup \B$ and let $u_k$ denote the number of permutations of length $k$ in $\U$. Then:
\[
\max (a_k , b_k) \leq u_k \leq a_k + b_k.
\]
It follows immediately that $\sw_{\U} = \max (\sw_{\A}, \sw_{\B})$. 

Using Theorems \ref{THM:Expectation} and \ref{THM:Concentration} it follows that:
\[
c_{\U} = \max ( c_{\A}, c_{\B} ).
\]
So $\U$ satisfies Conjecture \ref{conj01}. The argument for direct sum is almost equally simple based on the observation that the number of permutations of length $n$ that belong to the direct sum is at least equal to the maximum of $a_k b_{n-k}$ for $0 \leq k \leq n$ and at most $n$ times this quantity. This establishes that the Stanley-Wilf limit for the direct sum is the maximum $\sw_{\A}$ and $\sw_{\B}$. The same result follows for $c_{\A \oplus \B}$ using Theorems \ref{THM:Expectation} and \ref{THM:Concentration}.

Now let $\J = \juxt(\A, \B)$ and let $j_k$ denote the number of permutations of length $k$ in $\J$. Then:
\[
j_n \leq \sum_{k = 0}^n {n \choose k} a_k b_{n-k} \leq (n+1) j_n
\]
since each permutation in $\J$ of length $n$  has at least one and at most $n+1$ representations as a juxtaposition of a permutation in $\A$ and one in $\B$. Observing the similarity of the middle expression to the binomial expansion of $(\sw_{\A} + \sw_{\B})^n$ one obtains by taking $n$th roots that
\[
\sw_{\J} = \sw_{\A} + \sw_{\B}.
\]
On the other hand the arguments of Steele in \cite{Steele:DM81} applied to juxtapositions of arbitrary pattern avoidance classes rather than just the increasing and decreasing class together with the theorems above also yield:
\[
c_{\J} = c_{\A} + c_{\B}.
\]
So Conjecture \ref{conj01} holds for $\juxt(\A, \B)$.

Finally, suppose that $\A \cap \B$ is finite, and let $s$ be the maximum length of a permutation in $\A \cap \B$. Let $\M = \merge(\A, \B)$ and let $j_k$ denote the number of permutations of length $k$ in $\J$. Consider first $c_{\M}$. Given any permutation $\pi$, the merge of any of its longest $\A$-sequences with any of its longest $\B$-sequences is an $\M$-sequence in $\pi$. Conversely, any $\M$-sequence is a merge of an $\A$-sequence and a $\B$-sequence. So
\[
L_{\A}(\pi) + L_{\B}(\pi) - s \leq L_{\M}(\pi) \leq L_{\A}(\pi) + L_{\B}(\pi).
\]
Since $s$ is fixed, on taking expectation dividing by $2 \sqrt{n}$ and taking a limit as $n \rightarrow \infty$ we obtain 
\[
\sqrt{c_{\M}} = \sqrt{c_{\A}} + \sqrt{c_{\B}}.
\]
So $c_{\M} = (\sqrt{c_{\A}} + \sqrt{c_{\B}})^2$.

To form the merge of a permutation of length $k$ and one of length $n-k$ we choose $k$ positions for the first permutation and also $k$ values. So:
\[
m_n \leq \sum_{k=0}^n {n \choose k}^2 a_k b_{n-k}.
\]
If $\pi \in \M$ is represented in two different ways as the merge of a permutation in $\A$ and one in $\B$ then the number of positions in which the two merges disagree (in the sense that one is representing the value in that position as part of an $\A$-permutation and the other as part of a $\B$-permutation) is at most $2s$. This is because if we consider positions represented in the first merge as coming from an $\A$-permutation and in the second as coming from a $\B$-permutation, the pattern occurring at these positions belongs to $\A \cap \B$ (and likewise for the opposite case). So, labeling positions in a merge by type, the merges representing $\pi$ correspond to a set of diameter at most $2s$ in the Hamming metric on the cube $\{0,1\}^n$. Kleitman \cite{Kleitman:JCT66} proved that the maximum number of points in such a set is attained in a ball of radius $s$, that is $v(n,s) = \sum_{t = 0}^s {n \choose s}$. So
\[
m_n \leq \sum_{k=0}^n {n \choose k}^2 a_k b_{n-k} \leq v(n,s) m_n.
\]
Since $s$ is fixed, $v(n,s)^{1/n} \rightarrow 1$ as $n \rightarrow \infty$. Therefore
\[
\sw_{\M} = \lim_{n \rightarrow \infty} \left( \sum_{k=0}^n {n \choose k}^2 a_k b_{n-k} \right)^{1/n}.
\]
Considering the similarity of the square root of each term of the sum in the expression above to a term in the expansion of $(\sqrt{\sw_{\A}} + \sqrt{\sw_{\B}})^n$ we obtain:
\[
\sw_{\M} = \left( \sqrt{\sw_{\A}} + \sqrt{\sw_{\B}} \right)^2
\]
and therefore $\M$ satisfies Conjecture \ref{conj01}.
\end{proof}

The results of the preceding proposition apply in part to the weaker version of Conjecture \ref{conj01} which only asserts the equality of two limits superior. Namely, following the details of the proof we see that if $\A$ and $\B$ satisfy the weaker form then so do their union and direct sum, while if one satisfies the stronger form and the other the weaker form then their juxtaposition and merge (subject to finite intersection) also satisfy the weaker form. Note also that the condition that $\A \cap \B$ be finite is equivalent to asserting that for some $k$ and $k'$, $k(k-1) \cdots 321 \not \in \A \cap \B$ and $123 \cdots k' \not \in \A \cap \B$. The ``interesting'' cases are when neither $A$ nor $B$ is finite, say $\A \subseteq \avoid(k(k-1) \cdots 321)$ and $\B \subseteq \avoid(123 \cdots k')$.

Let $\A$ be any pattern avoidance class and define $\cycle (\A)$ to be the set of permutations obtained by taking all the cyclic rotations of elements of $\A$. It is easily verified that $\cycle (\A)$ is also a pattern avoidance class.

\begin{proposition}
If $\A$ satisfies Conjecture \ref{conj01}  then so does $\cycle (\A)$.
\end{proposition}

\begin{proof}
Let $\A$ be a pattern avoidance class satisfying Conjecture \ref{conj01} and let $\R = \cycle (A)$. Since $\R \supseteq \A$ and the number of elements of $\R$ of length $n$ is at most $n$ times the number of elements of $\A$ of length $n$, $\sw_{\R} = \sw_{\A}$.

However, the definition of $c_{\R}$ and Theorem \ref{THM:Concentration} also imply that $c_{\R} = c_{\A}$. To see this observe that the length of any $\LS{\R}$ sequence in a permutation $\pi$ is the maximum of the lengths of an $\LS{\A}$ sequence in the cyclic rotations of $\pi$. The $n!$ permutations in $\S_n$ are partitioned into $(n-1)!$ classes of $n$ elements under cyclic rotation, and the mean length of an $\LS{\R}$ sequence is the mean of the maximum lengths of the $\LS{\A}$ sequences in each of these equivalence classes. Thus the mean length of an $\LS{\R}$ sequence is not greater than $M_n$, the mean length of an $\LS{\A}$ sequence over the $(n-1)!$ permutations whose $\LS{A}$ sequences are longest. Choose $1/3 < \alpha < 1/2$ and $0 < \beta < 3 \alpha - 1$. Then, using Theorem \ref{THM:Concentration}:
\[
M_n \leq \left( \expectation L_{\A} (\Pi_n) + n^{\alpha} \right) + n^2 \exp (-n^\beta).
\]
The final term arises from the upper bound of $n$ for the length of the longest $\A$ subsequence. In particular:
\[
\lim_{n \to \infty} \frac{M_n}{\sqrt{n}} \leq 2 \sqrt{c_{\A}}.
\]
This gives $c_{\R} \leq c_{\A}$ and the reverse inequality is of course trivial. Thus $\R$ satisfies Conjecture \ref{conj01} as claimed.
\end{proof}

Again, the ``weak'' form of this proposition is valid.

\section{Experimental results}
\label{SEC:Experiment}

The results of the previous section provide some evidence in favour of Conjecture \ref{conj01}, or more conservatively, provide evidence that for a fairly wide set of pattern avoidance classes $\A$, $c_{\A} = \sw_{\A}$. However, we have no significant evidence in favour of this conjecture based on an exact computation of $c_{\A}$ for any classes other than those which can be produced from the classes $\avoid(k(k-1)\cdots321)$ or $\avoid(123\cdots(m-1)m)$ by the constructions described in the previous section.

In this section we consider two specific pattern avoidance classes:
\begin{itemize}
\item
The {\em layered} permutations, $\L = \avoid(231, 312)$, consisting of all permutations of the form $D_1 \oplus D_2 \oplus \cdots \oplus D_k$ where each of $D_1$ through $D_k$ is a descending permutation. The number of permutations of length $n$ in $\L$ is $2^{n-1}$ and so $\sw_{\L} = 2$.
\item
The subclass $\L(2) = \avoid(231, 312, 321)$ of $\L$ formed by requiring that each $D_i$ contain at most two elements. The number of permutations of length $n$ in $\L(2)$ is equal to the $n$th Fibonacci number, so $\sw_{\L(2)} = (1 + \sqrt{5})/2$. Kaiser and Klazar \cite{KK:EJC0203} proved that $\L(2)$ is the smallest pattern avoidance class whose Stanley-Wilf limit is strictly greater than $1$.
\end{itemize}

Results in \cite{AAADHHO:AJC03} give dynamic programming algorithms for solving the longest subsequence problem for both $\L$ and $\L(2)$ whose complexity is $O(n^2 \log n)$ where $n$ is the length of the input permutation. We have been able to improve the latter algorithm, based on a tableau style method to result in a complexity of $O(n \log n)$. Both algorithms were implemented and a long period random number generator was used to provide experimental data concerning the values $c_{\L}$ and $c_{\L(2)}$.

For $\L(2)$ we present data based on permutations of length $2^k \times 10^4$ for $0 \leq k \leq 7$. For each value of $k$, $1000$ random permutations of that length were generated and the length of the longest $\L(2)$ subsequences was computed. Table \ref{tableL2} shows the mean, sample standard deviation, and resulting estimates of $c_{\L(2)}$ based on these simulations. We would be forced to classify a person who believed in the truth of Conjecture \ref{conj01} based on this data for $\L(2)$ as an optimist. If the estimates are indeed converging to $\sw_{\L(2)}$ then they are not yet within $6\%$ of their final limit at $n = 128 \times 10^4$. By contrast, for this value of $n$ the estimate for $c_{\I}$ (whose actual value is $1$) is approximately $0.985$.

\begin{table}
\begin{tabular}{c|cc|c} 
Length & $\mu$ & $\sigma$ & $ \sim c_{\L(2)}$ \\ \hline
$1 \times 10^4$ & 239.3 & 4.5 & 1.431 \\
$2 \times 10^4$ & 340.7 & 5.2 & 1.451 \\
$4 \times 10^4$ & 484.7 & 6.1 & 1.468 \\
$8 \times 10^4$ & 688.4 & 6.4 & 1.481 \\
$16 \times 10^4$ & 978.1 & 7.1 & 1.495 \\
$32 \times 10^4$ & 1386.8 & 8.3 & 1.503 \\
$64 \times 10^4$ & 1965.3 & 9.3 & 1.510 \\
$128 \times 10^4$ & 2785.3 & 10.2 & 1.515 \\
\end{tabular}
\caption{Summary data for the mean, $\mu$, and standard deviation $\sigma$ of the length of the longest $\L(2)$ subsequences of random permutations together with corresponding estimates of $c_{\L(2)}$.}
\label{tableL2}
\end{table}

Because of the slower running time and increased space requirements required by the algorithm for finding longest layered subsequences data for $\L$ is based on permutations of length $2^k \times 10^2$ for $0 \leq k \leq 7$. As for $\L(2)$, $1000$ random permutations of each length were analysed and the results are presented in Table \ref{tableL}. The data for this class do not require as much optimism as the $\L(2)$ data to be viewed as support for Conjecture 1. 
\begin{table}
\begin{tabular}{c|cc|c} 
Length & $\mu$ & $\sigma$ & $ \sim c_{\L}$ \\ \hline
$1 \times 10^2$ & 23.8 & 1.8 & 1.418 \\
$2 \times 10^2$ & 34.8 & 2.2 & 1.517 \\
$4 \times 10^2$ & 50.6 & 2.5 & 1.602 \\
$8 \times 10^2$ & 73.4 & 3.0 & 1.682 \\
$16 \times 10^2$ & 105.2 & 3.3 & 1.730 \\
$32 \times 10^2$ & 150.7 & 4.0 & 1.774 \\
$64 \times 10^2$ & 215.9 & 4.4 & 1.821 \\
$128 \times 10^2$ & 307.5 & 4.9 & 1.847 \\
\end{tabular}
\caption{Summary data for the mean, $\mu$, and standard deviation $\sigma$ of the length of the longest $\L$ subsequences of random permutations together with corresponding estimates of $c_{\L}$.}
\label{tableL}
\end{table}

\section{Conclusions}

We have illustrated that, broadly speaking, the $L_{\A}$ statistic on $\S_n$ has the same general qualities as the more well known $\L_{\I}$ statistic. Certainly a great deal more work is required before one could begin to appreciate the former statistic at the level of detail that is known concerning the latter. Not surprisingly, it appears that there is a connection between $\L_{\A}$ and the growth rate of the pattern avoidance class $\A$. Conjecture \ref{conj01} proposes a precise form for this connection and the results of Section \ref{SEC:Theory} establish a relatively wide collection of classes for which the Conjecture \ref{conj01} holds. One extension to these results concerns classes $\A$ for which $\sw_{\A} = 1$. These are (essentially) classified in \cite{KK:EJC0203} and it is only a matter of somewhat tedious routine to confirm that in all these cases $c_{\A} = 1$ also.

However, the experimental results of Section \ref{SEC:Experiment} can be viewed as providing evidence against Conjecture \ref{conj01}. In particular the class $\L(2)$ seems a likely contender as a possible counterexample. The structure of the permutations in this class is sufficiently simple that there may be some chance of carrying out a more detailed analysis of $\L_{\L(2)}$ with a view to explicitly calculating, or providing bounds for, $c_{\L(2)}$.

Notable by its omission from either Section \ref{SEC:Theory} or Section \ref{SEC:Experiment} is the class $\avoid (312)$. This class has Stanley-Wilf limit $4$ (as does every class defined by avoiding a single three element pattern since all such classes are enumerated by the Catalan numbers). A polynomial time algorithm for the longest subsequence problem based on this class is given in \cite{AAADHHO:AJC03} but its complexity on permutations of length $n$ is $O(n^5)$ which makes it impractical for experiments of the size required to produce even vaguely convincing evidence. The goal of producing such evidence would seem to require finding, even on an {\em ad hoc} basis some collection of classes for which the longest subsequence problem can be solved algorithmically in reasonable time (basically, at worst quadratic) and/or developing better algorithms for classes such as $\avoid(312)$.

\bibliographystyle{plain}
\bibliography{LAS}

\end{document}